\documentclass{article}
\usepackage{graphicx}
\usepackage{amsmath}
\usepackage{amsfonts}
\usepackage{amssymb}
\newtheorem{theorem}{Theorem}

\newtheorem{corollary}[theorem]{Corollary}

\newtheorem{definition}[theorem]{Definition}

\newtheorem{lemma}[theorem]{Lemma}

\newtheorem{proposition}[theorem]{Proposition}

\newenvironment{proof}[1][Proof]{\textbf{#1.} }{\ \rule{0.5em}{0.5em}}

\begin{document}

\title{Intersection theory of coassociative submanifolds in $G_{2}$-manifolds and
Seiberg-Witten invariants}
\author{Naichung Conan Leung and Xiaowei Wang}
\maketitle
\begin{abstract}
We study the problem of counting instantons with coassociative boundary
condition in (almost) $G_{2}$-manifolds. This is analog to the open
Gromov-Witten theory for counting holomorphic curves with Lagrangian boundary
condition in Calabi-Yau manifolds. We explain its relationship with the
Seiberg-Witten invariants for coassociative submanifolds.
\end{abstract}

\bigskip

Intersection theory of Lagrangian submanifolds is an essential part of the
symplectic geometry. By counting the number of holomorphic disks bounding
intersecting Lagrangian submanifolds, Floer and others defined the celebrated
Floer homology theory. It plays an important role in mirror symmetry for
Calabi-Yau manifolds and string theory in physics. In M-theory, Calabi-Yau
threefolds are replaced by seven dimensional $G_{2}$-manifolds $M$ (i.e.
oriented Octonion manifolds \cite{Leung RG over A}). The analog of holomorphic
disks (resp. Lagrangian submanifolds) are instantons or associative
submanifolds (resp. coassociative submanifolds or branes) in $M$ \cite{Lee
Leung Instanton Brane}. An important project is to count the number of
instantons with coassociative boundary conditions. In particular we want to
study the following problem.

\bigskip

\textbf{Problem}: Given two nearby coassociative submanifolds $C$ and
$C^{\prime}$ in a (almost) $G_{2}$-manifold $M$. Relate the number of
instantons in $M$ bounding $C\cup C^{\prime}$ to the Seiberg-Witten invariants
of $C$.

\bigskip

The basic reason is a coassociative submanifold $C^{\prime}$ which is
infinitesimally close to $C$ corresponds to a symplectic form on $C$ which
degenerates along $C\cap C^{\prime}$. Instantons bounding $C\cup C^{\prime}$
would become holomorphic curves on $C$ modulo bubbling. By the work of Taubes,
we expect that the number of such instantons is given by the Seiberg-Witten
invariant of $C$.

In this paper we treat the special case when $C$ and $C^{\prime}$ are
disjoint, i.e. $C$ is a symplectic four manifold. Recall that Taubes showed
that the Seiberg-Witten invariants of such a $C$ is given by the Gromov-Witten
invariants \cite{Ta Gr Sw} of $C$. Our main result\ is following theorem.

\begin{theorem}
Suppose that $M$ is an (almost) $G_{2}$-manifold and $\left\{  C_{t}\right\}
$ is an one parameter family of coassociative submanifolds in $M$ such that
$dC_{t}/dt|_{t=0}$ is nonvanishing.

If $\left\{  A_{t}\right\}  $ is any one parameter family of instantons in $M$
satisfying
\[
\partial A_{t}\subset C_{t}\cup C_{0}\text{ \ and }\lim_{t\rightarrow0}%
A_{t}=\Sigma\text{ in }C^{1}\text{-topology,}%
\]
then $\Sigma$ is a $J$-holomorphic curve in $C_{0}$. \ 

Conversely, suppose that $\Sigma$ is a \emph{regular\ }$J$-holomorphic curve
in $C_{0}$,\ then it is the limit of a family of instantons $A_{t}$'s as
described above.
\end{theorem}

\bigskip

A few remarks are in order: First, counting such small instantons is basically
a problem in four manifold theory because of Bryant's result \cite{Bryant}
which says that the zero section $C$ in $\Lambda_{+}^{2}\left(  C\right)  $ is
always a coassociative submanifold for some incomplete $G_{2}$-metric on its
neighborhood provided that the bundle $\Lambda_{+}^{2}\left(  C\right)  $ is
topologically trivial. Second, when $C$ and $C^{\prime}$ are not disjoint, the
above theorem should still hold true. However using the present approach to
prove it would require a good understanding of the Seiberg-Witten theory on
any four manifold with a degenerated symplectic form as in Taubes program.
Third, when $C\,$\ and $C^{\prime}$ are not close to each other then we have
to take into account the bubbling phenomenon which has not been established
yet. Nevertheless, one would expect that if the volume of $A_{t}$'s are small,
then bubbling cannot occur, thus they would converge to a holomorphic curve in
$C_{0}$.

\section{Review of Symplectic Geometry}

Given any symplectic manifold $\left(  X,\omega\right)  $ of dimension $2n$,
there exists a compatible metric $g$ so that the equation
\[
\omega\left(  u,v\right)  =g\left(  Ju,v\right)
\]
defines a Hermitian almost complex structure
\[
J:T_{X}\rightarrow T_{X}\text{,}%
\]
that is $J^{2}=-id$ and $g\left(  Ju,Jv\right)  =g\left(  u,v\right)  $.

A \emph{holomorphic curve}, or \emph{instanton}, is a two dimensional
submanifold $\Sigma$ in $X$ whose tangent bundle is preserved by $J$.
Equivalently $\Sigma$ is calibrated by $\omega,$ i.e. $\omega|_{\Sigma
}=vol_{\Sigma}$. By \emph{counting} the number of instantons in $X$, one can
define a highly nontrivial invariant for the symplectic structure on $X$,
called the Gromov-Witten invariant.

When the instanton $\Sigma$ has nontrivial boundary, then the corresponding
free boundary value problem would require $\partial\Sigma$ to lie on a
\emph{Lagrangian submanifold} $L$ in $X$, i.e. $\dim L=n$ and $\omega|_{L}=0
$. Floer studied the intersection theory of Lagrangian submanifolds and
defined the Floer homology group $HF\left(  L,L^{\prime}\right)  $ under
certain assumptions.

Suppose that $X$ is a Calabi-Yau manifold, i.e. the holonomy group of the
Levi-Civita connection is inside $SU\left(  n\right)  $, equivalently $J$ is
an integrable complex structure on $X$ and there exists a holomorphic volume
form $\Omega_{X}\in\Omega^{n,0}\left(  X\right)  $ on $X$ satisfying
$\Omega_{X}\bar{\Omega}_{X}=C_{n}\omega^{n}$. Under the mirror symmetry
transformation, $HF\left(  L,L^{\prime}\right)  $ is expected to correspond to
the Dolbeault cohomology group of coherent sheaves in the mirror Calabi-Yau manifold.

A Lagrangian submanifold $L$ in $X$ is called a \emph{special Lagrangian
submanifold }with phase zero (resp. $\pi/2$) if $\operatorname{Im}\Omega
_{X}|_{L}=0 $ (resp. $\operatorname{Re}\Omega_{X}|_{L}=0$). Such a $L$ is
calibrated by $\operatorname{Re}\Omega_{X}|_{L}$ (resp. $\operatorname{Im}%
\Omega_{X}|_{L}$). They play important roles in the Strominger-Yau-Zaslow
mirror conjecture for Calabi-Yau manifolds \cite{SYZ}.

When $X$ is a Calabi-Yau threefold, there are conjectures of Vafa and others
(e.g. \cite{Mina Vafa}\cite{GV2}) that relates the (partially defined) open
Gromov-Witten invariant of the number of instantons with Lagrangian boundary
condition to the large $N$ Chern-Simons invariants of knots in three manifolds.

\section{Counting Instantons in (almost) $G_{2}$-manifolds}

Notice that a real linear homomorphism $J:\mathbb{R}^{m}\rightarrow
\mathbb{R}^{m}$ being a Hermitian complex structure on $\mathbb{R}^{m}$ is
equivalent to the following conditions: for any vector $v\in\mathbb{R}^{m}$ we
have (i) $Jv$ is perpendicular to both $v$ and (ii) $\left|  Jv\right|
=\left|  v\right|  $. We can generalize $J$ to involve more than one vector.
We call a skew symmetric homomorphism
\[
\times:\mathbb{R}^{m}\otimes\mathbb{R}^{m}\rightarrow\mathbb{R}^{m}%
\]
a (2-fold) \emph{vector cross product} if it satisfies
\[%
\begin{array}
[c]{cl}%
\text{(i) } & \left(  u\times v\right)  \,\text{is perpendicular to both
}u\text{ and }v\text{, and}\\
\text{(ii) } & \left|  u\times v\right|  =\text{Area of parallelogram spanned
by }u\text{ and }v\text{.}%
\end{array}
\]
The obvious example of this is the standard vector product on $\mathbb{R}^{3}
$. By identifying $\mathbb{R}^{3}$ with $\operatorname{Im}\mathbb{H}$, the
imaginary part of the quaternion numbers, we have
\[
u\times v=\operatorname{Im}u\bar{v}.
\]
The same formula defines a vector cross product on $\mathbb{R}^{7}%
=\operatorname{Im}\mathbb{O}$, the imaginary part of the octonion numbers.
Brown and Gray \cite{Gray VectorCrossProd} showed that these two are the only
possible vector cross product structures on $\mathbb{R}^{m}$ up to isomorphisms.

Suppose that $M$ is a seven dimensional Riemannian manifold with a vector
cross product $\times$ on each of its tangent spaces. The analog of the
symplectic form is a degree three differential form $\Omega$ on $M$ defined as
follow:
\[
\Omega\left(  u,v,w\right)  =g\left(  u\times v,w\right)  \text{.}%
\]

\begin{definition}
Suppose that $\left(  M,g\right)  $ is a Riemannian manifold of dimension
seven with a vector cross product $\times$ on its tangent bundle. Then (1) $M$
is called an \textbf{almost }$G_{2}$\textbf{-manifold} if $d\Omega=0$ and (2)
$M$ is called a\textbf{\ }$G_{2}$\textbf{-manifold} if $\nabla\Omega=0$.
\end{definition}

It can be proven that the condition $\nabla\Omega=0$ is equivalent to $\Omega$
being a harmonic form, i.e. $\Delta\Omega=0$. Furthermore $M$ is a $G_{2}%
$-manifold if and only if its holonomy group is inside the exceptional Lie
group $G_{2}=Aut\left(  \mathbb{O}\right)  $. The geometry of $G_{2}%
$-manifolds can be interpreted as the symplectic geometry on its knot space
(see e.g. \cite{Lee Leung Instanton Brane}, \cite{Movshev}).

For example, if $\left(  X,\omega_{X}\right)  $ is a Calabi-Yau threefold with
a holomorphic volume form $\Omega_{X}$, then the product manifold $M=X\times
S^{1}$ is a $G_{2}$-manifold with
\[
\Omega=\operatorname{Re}\Omega_{X}+\omega_{X}\wedge d\theta\text{.}%
\]

Next we define the analogs of holomorphic curves and Lagrangian submanifolds
in the $G_{2}$ setting.

\begin{definition}
Suppose that $A$ is a three dimensional submanifold of an almost $G_{2}%
$-manifold $M$. We call $A$ an \textbf{instanton }or \textbf{associative
submanifold}, if $A$ is preserved by the vector cross product $\times$.
\end{definition}

Harvey and Lawson \cite{Harvey Lawson} showed that $A\subset M$ is an
instanton if and only if $A$ is calibrated by $\Omega$, i.e. $\Omega
|_{A}=vol_{A}$.

In M-theory, associative submanifolds are also called \emph{M2-branes}. For
example when $M=X\times S^{1}$ with $X$ a Calabi-Yau threefold, $\Sigma\times
S^{1}$ (resp. $L\times\left\{  p\right\}  $) is an instanton in $M $ if and
only if $\Sigma$ (resp. $L$) is a holomorphic curve (resp. special Lagrangian
submanifold with zero phase) in $X$.

A natural interesting question is to count the number of instantons in $M$. In
the special case of $M=X\times S^{1}$, these numbers are reduced to the
conjectural invariants proposed by Joyce \cite{Joyce Count SLag} by counting
special Lagrangian submanifolds in Calabi-Yau threefolds. This problem has
been discussed by many physicists. For example Harvey and Moore discussed in
\cite{Harvey Moore} the mirror symmetry aspects of these invariants; Aganagic
and Vafa in \cite{Mina Vafa} related these invariants to the open
Gromov-Witten invariants for local Calabi-Yau threefolds; Beasley and Witten
argued in \cite{BeasleyWitten} that when there is a moduli of instantons, then
one should count them using the Euler characteristic of the moduli space. In
this paper we count the number of instantons with boundary lying on a
coassociative submanifold in $M$. The compactness issues of the moduli of
instantons is a very challenging problem because the dimension of an instanton
is bigger than two. This makes it very difficult to define an honest invariant
by counting instantons.

\bigskip\ 

When an instanton $A$ has a nontrivial boundary, $\partial A\neq\phi,$ one
should require it to lie inside a \emph{brane }or a \emph{coassociative
submanifold} \cite{Lee Leung Instanton Brane}, i.e. submanifolds in $M$ where
the restriction of $\Omega$ is zero and have the largest possible dimension.
Branes are the analog of Lagrangian submanifolds in symplectic geometry.

\begin{definition}
Suppose that $C$ is a four dimensional submanifold of an almost $G_{2}%
$-manifold $M$. We call $C$ a \textbf{coassociative submanifold} if
\[
\Omega|_{C}=0\text{ and }\dim C=4.
\]
\end{definition}

For example when $M=X\times S^{1}$ with $X$ a Calabi-Yau threefold, $H\times
S^{1}$ (resp. $C\times\left\{  p\right\}  $) is a coassociative submanifold in
$M$ if and only if $H$ (resp. $C$) is a special Lagrangian submanifold with
phase $\pi/2$ (resp. complex surface) in $X$. In \cite{Lee Leung Instanton
Brane} J.H. Lee and the first author showed that the isotropic knot space
$\mathcal{\hat{K}}_{S^{1}}X$ of $X$ admits a natural holomorphic symplectic
structure. Moreover $\mathcal{\hat{K}}_{S^{1}}H$ (resp. $\mathcal{\hat{K}%
}_{S^{1}}C$) is a complex Lagrangian submanifold in $\mathcal{\hat{K}}_{S^{1}%
}X$ with respect to the complex structure $J$ (resp. $K$).

Constructing special Lagrangian submanifolds with zero phase in $X$ with
boundaries lying on $H$ (resp. $C$) corresponds to the Dirichlet (resp.
Neumann) free boundary value problem for minimizing volume among Lagrangian
submanifolds as studied by Schoen and others. For a general $G_{2}$-manifold
$M$, the natural free boundary value for an instanton is a coassociative
submanifold. Similar to the intersection theory of Lagrangian submanifolds in
symplectic manifolds. We propose to study the following problem: Count the
number of instantons in $G_{2}$-manifolds bounding two coassociative submanifolds.

The product of a coassociative submanifold with a two dimensional plane inside
the eleven dimension spacetime $M\times\mathbb{R}^{3,1}$ is called a
\emph{D5-brane }in M-theory. Counting the number of M2-branes between two
D5-branes has also been studied in the physics literatures.

In general this is a very difficult problem. For instance, counting $S^{1}%
$-invariant instantons in $M=X\times S^{1}$ is the open Gromov-Witten
invariants. However when the two coassociative submanifolds $C$ and
$C^{\prime}$ are close to each other, we can relate the number of instantons
between them to the Seiberg-Witten invariant of $C$.

\section{Relationships to Seiberg-Witten invariants}

To determine the number of instantons between nearby coassociative
submanifolds, we first recall the deformation theory of coassociative
submanifolds $C$ inside any $G_{2}$-manifold $M$, as developed by McLean
\cite{McLean}. Given any normal vector $v\in N_{C/M}$, the interior product
$\iota_{v}\Omega$ is naturally a self-dual two form on $C$ because of
$\Omega|_{C}=0$. This gives a natural identification,
\begin{align*}
&  N_{C/M}\overset{\simeq}{\rightarrow}\Lambda_{+}^{2}\left(  C\right) \\
v &  \rightarrow\eta_{0}=\iota_{v}\Omega\text{.}%
\end{align*}
Furthermore infinitesimal deformations of coassociative submanifolds are
parametrized by self-dual harmonic two forms $\eta_{0}\in H_{+}^{2}\left(
C\right)  $, and they are always unobstructed. Notice that the zero set of
$\eta_{0}$ is the intersection of $C$ with a infinitesimally near
coassociative submanifold, that is
\[
\left\{  \eta_{0}=0\right\}  =\lim_{t\rightarrow0}\left(  C\cap C_{t}\right)
\text{,}%
\]
where $C=C_{0}$ and $\eta_{0}=dC_{t}/dt|_{t=0}$.

Since
\[
\eta_{0}\wedge\eta_{0}=\eta_{0}\wedge\ast\eta_{0}=\left|  \eta_{0}\right|
^{2}\ast1\text{,}%
\]
$\eta_{0}$ defines a natural symplectic structure on $C^{reg}:=C\backslash
\left\{  \eta_{0}=0\right\}  $.\ If we normalize $\eta_{0}$,
\[
\eta=\eta_{0}/\left|  \eta_{0}\right|  \text{,}%
\]
then the equation
\[
\eta\left(  u,v\right)  =g\left(  Ju,v\right)
\]
defines a Hermitian almost complex structure on $C^{reg}$.

The next lemma says that when two coassociative submanifolds $C$ and
$C^{\prime}$ come together, then the limit of instantons bounding them will be
a holomorphic curve $\Sigma$ in $C^{reg}$ with boundary $C\cap C^{\prime}$.

\begin{proposition}
Suppose that $C_{t}$ is an one parameter family of coassociative submanifolds
in a $G_{2}$-manifold $M$. Suppose that $A_{t}$ is a family of instantons in
$M$ bounding $C_{0}\cup C_{t}$ for nonzero $t$ and
\[
\lim_{t\rightarrow0}A_{t}=\Sigma
\]
exists in $C^{1}$-topology. Then $\Sigma$ is a $J$-holomorphic curve in
$C_{0}^{reg}$ with boundary $C_{0}\backslash C_{0}^{reg}$.
\end{proposition}

\begin{proof}
For simplicity we assume that $\eta_{0}=dC_{t}/dt|_{t=0}$ is nowhere
vanishing. Let us denote the boundary component of $A_{t}$ in $C_{0}$ as
$\Sigma_{t}$ and the unit normal vector field for $\Sigma_{t}$ in $A_{t}$ as
$n_{t}$. Note that $n_{t}$ is perpendicular to $C_{0}$. This is because $A_{t}
$ being preserved by the vector cross product implies that
\[
n_{t}=u\times v,
\]
for some tangent vectors $u$ and $v$ in $\Sigma_{t}$, therefore given any
tangent vector $w$ along $C_{0}$, we have
\[
g\left(  n_{t},w\right)  =g\left(  u\times v,w\right)  =\Omega\left(
u,v,w\right)  =0.
\]
The last equality follows from $C_{0}$ being coassociative and $\Sigma
_{t}\subset C_{0}$. Using this and the fact that $A_{t}$ bounds $C_{0}\cup
C_{t}$ with $\lim_{t\rightarrow0}C_{t}=C_{0}$, i.e. $n_{t}$ is \emph{pointing
towards} $C_{t}$, we obtain%

\[
\lim_{t\rightarrow0}n_{t}=\eta|_{\Sigma}\text{.}%
\]
Therefore $\Sigma=\lim_{t\rightarrow0}\Sigma_{t}$ is a holomorphic curve in
$C_{0}$ with respect to the almost complex structure $J$ defined by
$\eta\left(  u,v\right)  =g\left(  Ju,v\right)  $.
\end{proof}

\bigskip

The reverse of the above proposition should also hold true. The Lagrangian
analog of it is proven by Fukaya and Oh in \cite{Fukaya Oh}. On the other
hand, by the celebrated work of Taubes, we expect that the number of such open
holomorphic curves in $C_{0}$ equals to the Seiberg-Witten invariant of
$C_{0}$. We conjecture the following statement.

\bigskip

\textbf{Conjecture}: Suppose that $C$ and $C^{\prime}$ are nearby
coassociative submanifolds in a $G_{2}$-manifold $M$. Then the number of
instantons in $M$ with small volume and with boundary lying on $C\cup
C^{\prime}$ is given by the Seiberg-Witten invariants of $C$.

\bigskip

In the next section we will discuss the case when $C$ and $C^{\prime}$ do not
intersect. The basic ideas are (i) the limit of such instantons is a
holomorphic curve with respect to the (degenerated) symplectic form $\eta$ on
$C$ coming from its deformations as coassociative submanifolds and this
process can be reversed; (ii) the number of holomorphic curves in the four
manifold $C$ should be related to the Seiberg-Witten invariant of $C$ by the
work of Taubes (\cite{Ta ICM1998}, \cite{Ta SW GW deg}).

Suppose that $\eta$ is a self-dual two form on $C$ with constant length
$\sqrt{2}$, in particular it is a (non-degenerate) symplectic form, and
$\Sigma$ is a holomorphic curve in $C$, possibly disconnected. If $\Sigma$ is
\emph{regular} in the sense that the linearized operator $\bar{\partial}$
\ has trivial cokernel \cite{Ta Gr Sw}, then Taubes showed that the perturbed
Seiberg-Witten equations,
\begin{align*}
F_{a}^{+} &  =\tau\left(  \psi\otimes\psi^{\ast}\right)  -r\sqrt{-1}\eta,\\
D_{A\left(  a\right)  }\psi &  =0,
\end{align*}
have solutions for all sufficient large $r$. Here $a$ is a connection on the
complex line bundle $E$ over $C$ whose first Chern class equals the
Poincar\'{e} dual of $\Sigma$, $PD\left[  \Sigma\right]  $, $\psi$ is a
section of the twisted spinor bundle $S_{+}=E\oplus\left(  K^{-1}\otimes
E\right)  $ and $D_{A\left(  a\right)  }$ is the twisted Dirac operator. The
number of such solutions is the Seiberg-Witten invariant $SW_{C}\left(
\Sigma\right)  $ of $C$. Furthermore the converse is also true, thus Taubes
established an equivalence between Seiberg-Witten theory and Gromov-Witten
theory for symplectic four manifolds. This result has far reaching
applications in four dimensional symplectic geometry.

For a general four manifold $C$ with nonzero $b^{+}\left(  C\right)  $, using
a generic metric, any self-dual two form $\eta$ on $C$ defines a degenerate
symplectic form on $C$, i.e. $\eta$ is a symplectic form on the complement of
$\left\{  \eta=0\right\}  $, which is a finite union of circles (see
\cite{Kirby}\cite{Honda}). Therefore, one might expect to have a relationship
between the Seiberg-Witten of $C$ and the number of holomorphic curves with
boundaries $\left\{  \eta=0\right\}  $ in $C$. Part of this Taubes' program
has been verified in \cite{Ta ICM1998}, \cite{Ta SW GW deg}.

\section{Proof of the main theorem}

Suppose that $\eta$ is a nowhere vanishing self-dual harmonic two form on a
coassociative submanifold $C$ in a $G_{2}$-manifold $M$. For any holomorphic
curve $\Sigma$ in $C$, we want to construct an instanton in $M$ bounding $C$
and $C^{\prime}$, where $C^{\prime}$ is a small deformation of the
coassociative submanifold $C$ along the normal direction $\eta$. Notice that
$C$ and $C^{\prime}$ do not intersect. We will construct such an instanton
using a perturbation argument which requires a lower bound on the first
eigenvalue for the appropriate elliptic operator. Recall that the deformation
of an instanton is governed by a twisted Dirac operator. We will reinterpret
it as a complexified version of the Cauchy-Riemann operator.

\subsection{Deformation of instantons}

To construct an instanton $A$ in $M$ from a holomorphic curve $\Sigma$ in $C
$, we need to perturb an almost instanton $A^{\prime}$ to a honest one using a
quantitative version of the implicit function theorem. Let us first recall the
deformation theory of instantons $A$ (\cite{Harvey Lawson} and \cite{Lee Leung
Instanton Brane}) in a Riemannian manifold $M$ with a parallel (or closed)
$r$-fold vector cross product
\[
\times:\Lambda^{r}T_{M}\rightarrow T_{M}\text{.}%
\]
In our situation, we have $r=2$. By taking the wedge product with $T_{M}$ we
obtain a homomorphism $\tau$,
\[
\tau:\Lambda^{r+1}T_{M}\rightarrow\Lambda^{2}T_{M}\cong\Lambda^{2}T_{M}^{\ast
}\text{,}%
\]
where the last isomorphism is induced from the Riemannian metric. As a matter
of fact, the image of $\tau$ lies inside the subbundle $\mathbf{g}_{M}^{\bot}$
which is the orthogonal complement of $\mathbf{g}_{M}\subset\mathbf{so}\left(
T_{M}\right)  \cong\Lambda^{2}T_{M}^{\ast}$, the bundle infinitesimal
isometries of $T_{M}$ preserving $\times$. That is,
\[
\tau\in\Omega^{r+1}\left(  M,\mathbf{g}_{M}^{\bot}\right)  .
\]

\begin{lemma}
(\cite{Harvey Lawson}, \cite{Lee Leung Instanton Brane}) An $r+1$ dimensional
submanifold $A\subset M$ is an instanton, i.e. preserved by $\times$, if and
only if
\[
\tau|_{A}=0\in\Omega^{r+1}\left(  A,\mathbf{g}_{M}^{\bot}\right)  .
\]
\end{lemma}

This lemma is important in describing deformations of an instanton. Namely it
shows that the normal bundle to an instanton $A$ is a twisted spinor bundle
over $A$ and infinitesimal deformations of $A$ are parametrized by twisted
harmonic spinors.

In our present situation, $M$ is a $G_{2}$-manifold. Using the interior
product with $\Omega$, we can identify $\mathbf{g}_{M}^{\bot}$ with the
tangent bundle $T_{M}$ and we can also characterize $\tau\in\Omega^{3}\left(
M,T_{M}\right)  $ by the following formula,
\[
\left(  \ast\Omega\right)  \left(  u,v,w,z\right)  =g\left(  \tau\left(
u,v,w\right)  ,z\right)  \text{.}%
\]
Therefore $A\subset M$ is an instanton if and only if $\ast\left(  \tau
|_{A}\right)  =0\in T_{M}|_{A}$. As a matter of fact, if $A$ is already close
to be an instanton, then we only need the normal components of $\ast\left(
\tau|_{A}\right)  $ to vanish.

\begin{proposition}
\label{1Prop alm instanton}There is a positive constant $\delta$ such that for
any three dimensional linear subspace $A$ in $M\cong\operatorname{Im}%
\mathbb{O}$ with $\left|  \tau|_{A}\right|  <\delta$, $A$ is an instanton if
and only if $\ast\left(  \tau|_{A}\right)  \in T_{A}$.
\end{proposition}

\begin{proof}
McLean \cite{McLean} observed that if $A_{t}$ is a family of linear subspaces
in $M\cong\mathbb{R}^{7}$ with $A_{0}$ an instanton, then
\[
\ast\left(  \frac{d\tau}{dt}|_{A_{t}}\right)  |_{t=0}\in N_{A_{0}/M}\subset
T_{M}|_{A_{0}}\text{.}%
\]
Explicitly, if we denote the standard base for $\mathbb{R}^{7}$ as $e_{i}$'s,
e.g. $e_{1}\times e_{2}=e_{3}$, then we can assume that $A$ is spanned by
$e_{1},e_{2}$ and $\tilde{e}_{3}=e_{3}+\sum_{i=4}^{7}t_{i}e_{i}$ for some
small $t_{i}$'s because the natural action of $G_{2}$ on the Grassmannian
$Gr\left(  2,7\right)  $ is transitive. Then an easy computation shows that
the normal component of $\ast\left(  \tau|_{A}\right)  $ in $N_{A/M}$ is given
by
\[
\ast\left(  \tau|_{A}\right)  ^{\perp}=-t_{5}\left(  e_{4}\right)  ^{\perp
}+t_{4}\left(  e_{5}\right)  ^{\perp}+t_{7}\left(  e_{6}\right)  ^{\perp
}-t_{6}\left(  e_{7}\right)  ^{\perp}.
\]
When $t_{j}$'s are all zero, we have $\left(  e_{j}\right)  ^{\perp}=e_{j}$
for $4\leq j\leq7$. In particular, they are linearly independent when $t_{j}
$'s are small. In that case, $\ast\left(  \tau|_{A}\right)  ^{\perp}=0 $ will
actually imply that $t_{j}=0$ for all $j$, i.e. $A$ is an instanton in $M$.
Hence the proposition.
\end{proof}

\bigskip

This proposition will be needed later when we perturb an almost instanton to
an honest one. We also need to identify the normal bundle $N_{A/M}$ to an
instanton $A$ with a twisted spinor bundle over $A$ as follow \cite{McLean}:
We denote $P$ the $SO\left(  4\right)  $-frame bundle of $N_{A/M}$. Using the
identification
\[
SO\left(  4\right)  =Sp\left(  1\right)  Sp\left(  1\right)  \rightarrow
SO\left(  \mathbb{H}\right)  ,
\]%
\[
\left(  p,q\right)  \cdot y=py\bar{q}\text{,}%
\]
the tangent bundle to $A$ can be identified as an associated bundle to $P$ for
the representation $SO\left(  4\right)  \rightarrow SO\left(
\operatorname{Im}H\right)  $, $\left(  p,q\right)  \cdot y=qy\bar{q}$. As a
result the spinor bundle $\mathbb{S}$ of $A$ is associated to the
representation $SO\left(  4\right)  \rightarrow SO\left(  \mathbb{H}\right)  $
given by $\left(  p,q\right)  \cdot y=y\bar{q}$. Hence we obtain
\[
N_{A/M}\cong\mathbb{S}\otimes_{\mathbb{H}}E\text{,}%
\]
where $E$ is the associated bundle to $P$ for the representation $SO\left(
4\right)  \rightarrow SO\left(  \mathbb{H}\right)  $ given by $\left(
p,q\right)  \cdot y=py$.

\subsection{Complexified Cauchy-Riemann equation}

Recall that the normal bundle to any instanton $A$ is a twisted spinor bundle
$\mathbb{S\otimes}_{\mathbb{H}}E$, or simply $\mathbb{S}$, over $A$. Let
$\mathcal{D}$ be the Dirac operator on $A$. If $V:=V^{a}\frac{\partial
}{\partial\omega^{a}}$ is a normal vector field to $A$ and we write the
covariant differentiation of $V$ as $\nabla\left(  V\right)  :=V_{i}^{a}%
\frac{\partial}{\partial\omega^{a}}\otimes\omega^{i}$, then by viewing $V$ as
a twisted spinor or a quaternion valued function on $A$,
\[
V=V^{4}+\mathbf{i}V^{5}+\mathbf{j}V^{6}+\mathbf{k}V^{7},
\]
we have,
\begin{align*}
\mathcal{D}V  &  =-\left(  V_{1}^{5}+V_{2}^{6}+V_{3}^{7}\right)
+\mathbf{i}\left(  V_{1}^{4}+V_{3}^{6}-V_{2}^{7}\right) \\
&  +\mathbf{j}\left(  V_{2}^{4}-V_{3}^{5}+V_{1}^{7}\right)  +\mathbf{k}\left(
V_{3}^{4}+V_{2}^{5}-V_{1}^{6}\right)  ,
\end{align*}
where $\mathcal{D}:\mathcal{=}\nabla_{1}\mathbf{i}+\nabla_{2}\mathbf{j}%
+\nabla_{3}\mathbf{k}$.

Let us first consider \ a simplified model, suppose that $A$ is a product
Riemannian three manifold $\left[  0,\varepsilon\right]  \times\Sigma$ with
coordinates $\left(  x_{1},z\right)  $ where $z=x_{2}+ix_{3}$. Let $e_{1}$ be
the unit tangent vector field on $A$ normal to $\Sigma$, namely along the
$x_{1}$-direction. We have
\[
\mathcal{D}=e_{1}\cdot\frac{\partial}{\partial x_{1}}+\bar{\partial},
\]
where $\bar{\partial}$ is the Dolbeault operator on the Riemann surface
$\Sigma$.

The Clifford multiplication of $e_{1}$ on $\mathbb{S}$ satisfies $e_{1}%
^{2}=-1$ and therefore we have an eigenspace decomposition $\mathbb{S}%
:=\mathbb{S}^{+}\oplus\mathbb{S}^{-}$ corresponding to eigenvalues $\pm i$.

If we write $V=\left(  u,v\right)  $ with $u=V^{4}+\mathbf{i}V^{5}%
\in\mathbb{S}^{+}$ and $v=V^{6}+\mathbf{i}V^{7}\in\mathbb{S}^{-}$, then we
have
\begin{align*}
\mathcal{D}V &  =\left(  \frac{\partial u}{\partial x_{1}}\mathbf{i}%
-\partial_{z}v\right)  +\left(  -\frac{\partial v}{\partial x_{1}}%
\mathbf{i+}\bar{\partial}_{z}u\right)  \cdot\mathbf{j}\\
&  =\left(  \left(  \frac{\partial u}{\partial x_{1}}+\mathbf{i}\partial
_{z}v\right)  +\left(  \frac{\partial v}{\partial x_{1}}+\mathbf{i}%
\bar{\partial}_{z}u\right)  \cdot\mathbf{j}\right)  \cdot\mathbf{i}\\
&  =\left[
\begin{array}
[c]{cc}%
\mathbf{i} & 0\\
0 & -\mathbf{i}%
\end{array}
\right]  \left(  \frac{\partial}{\partial x_{1}}+\left[
\begin{array}
[c]{cc}%
0 & \mathbf{i}\partial_{z}\\
\mathbf{i}\bar{\partial}_{z} & 0
\end{array}
\right]  \right)  \left[
\begin{array}
[c]{c}%
u\\
v
\end{array}
\right]  ,
\end{align*}
where \
\[
\bar{\partial}_{z}:=\frac{\partial}{\partial x_{2}}+\mathbf{i}\frac{\partial
}{\partial x_{3}}\text{ and }\partial_{z}:=\frac{\partial}{\partial x_{2}%
}-\mathbf{i}\frac{\partial}{\partial x_{3}}.
\]
We will also denote $\mathbf{i}\partial_{z}$ and $\mathbf{i}\bar{\partial}%
_{z}$ by $\partial^{+}$ and $\partial^{-}$ respectively. They are Dirac
operators on $\Sigma$ and they satisfy
\[
\partial^{+}=\left(  \partial^{-}\right)  ^{\ast}\text{.}%
\]

This implies that the Dirac equation $\mathcal{D}V=0$ \ is equivalent to the
following complexified Cauchy-Riemann equations,
\begin{align*}
\bar{\partial}_{z}u  &  =\frac{\partial v}{\partial x_{1}}\mathbf{i}\text{,}\\
\partial_{z}v  &  =\frac{\partial u}{\partial x_{1}}\mathbf{i.}%
\end{align*}

\subsection{\label{1Sec e v estimate}Eigenvalue estimates}

In this subsection we give a quantitative estimate of the eigenvalue of the
linearized operator for the above simplified model. \ To do so, we first
introduce the following function spaces for spinors $V=\left(  u,v\right)  $
over a product three manifold $A_{\varepsilon}=\left[  0,\varepsilon\right]
\times\Sigma$.

\begin{definition}
\item  Suppose that $\mathbb{S}$ is the spinor bundle over a product three
manifold $A_{\varepsilon}=\left[  0,\varepsilon\right]  \times\Sigma$. We
define function spaces
\[
H_{-}^{m}\left(  A_{\varepsilon}\right)  :=\left\{  V\in H^{m}\left(
A_{\varepsilon},\mathbb{S}\right)  |\text{ }v|_{\left\{  0\right\}
\times\Sigma}=v|_{\left\{  \varepsilon\right\}  \times\Sigma}=0\right\}
\]
and
\[
H_{+}^{m}\left(  A_{\varepsilon}\right)  :=\left\{  V\in H^{m}\left(
A_{\varepsilon},\mathbb{S}\right)  |\text{ }u|_{\left\{  0\right\}
\times\Sigma}=u|_{\left\{  \varepsilon\right\}  \times\Sigma}=0\right\}
\]
where $H^{m}\left(  A,\mathbb{S}\right)  $ is $m$-th Sobolev space of sections
of $\mathbb{S}.$
\end{definition}

It is well known (see for instance \cite{BW} Theorem 21.5) that the Dirac
operators
\[
\mathcal{D}_{\pm}:=\mathcal{D}|_{H_{\pm}^{1}}:H_{\pm}^{1}\left(
A_{\varepsilon},\mathbb{S}\right)  \rightarrow L^{2}\left(  A_{\varepsilon
},\mathbb{S}\right)
\]
give well-defined local elliptic boundary problems and the formal adjoint
$\mathcal{D}_{+}^{\ast}=\mathcal{D}_{-}$. We are going to obtain an estimate
for its first eigenvalue.

\begin{theorem}
\label{1Thm Eigen}\label{eigenvalue}Suppose $\lambda_{\partial^{+}}$ is the
first eigenvalue of \ $\Delta_{\Sigma}=\partial^{-}\partial^{+}$ in the space
$H^{1}\left(  \Sigma,\mathbb{S}^{+}\right)  $ and let
\[
\lambda_{\mathcal{D}}:=\inf_{V\in H_{-}^{1}}\frac{\int_{X}\left\Vert
\mathcal{D}V\right\Vert ^{2}}{\int_{X}\left\Vert V\right\Vert ^{2}}.
\]
Then
\[
\lambda_{\partial^{+}}\geq\lambda_{\mathcal{D}}\geq\min\left\{  \lambda
_{\partial^{+}},2/\varepsilon^{2}\right\}  .
\]
In particular, we have $\lambda_{\mathcal{D}}=\lambda_{\partial^{+}}$ for
small $\varepsilon$.
\end{theorem}

\begin{proof}
For $\forall$ $V=\left(  u,v\right)  \in H_{-}^{1}\left(  A_{\varepsilon
}\right)  \mathcal{\ }$ we have
\[
\left\langle \mathcal{D}V,\mathcal{D}V\right\rangle _{L^{2}}=\int_{\left[
0,\varepsilon\right]  \times\Sigma}\left\vert \frac{\partial V}{\partial
x_{1}}\right\vert ^{2}+2\left\langle \frac{\partial V}{\partial x_{1}%
},\left[
\begin{array}
[c]{cc}%
0 & \partial^{+}\\
\partial^{-} & 0
\end{array}
\right]  V\right\rangle +\left\Vert \partial^{-}v\right\Vert ^{2}+\left\Vert
\partial^{+}u\right\Vert ^{2}.
\]
Using the formula $\partial^{-}=\left(  \partial^{+}\right)  ^{\ast}$, we
have
\begin{align*}
&  \int_{\left[  0,\varepsilon\right]  \times\Sigma}\left\langle
\frac{\partial V}{\partial x_{1}},\left[
\begin{array}
[c]{cc}%
0 & \partial^{+}\\
\partial^{-} & 0
\end{array}
\right]  V\right\rangle \\
&  =\int_{\left[  0,\varepsilon\right]  \times\Sigma}\left\langle
\frac{\partial u}{\partial x_{1}},\partial^{+}v\right\rangle +\left\langle
\frac{\partial v}{\partial x_{1}},\partial^{-}u\right\rangle \\
&  =\int_{\left[  0,\varepsilon\right]  \times\Sigma}\left\langle \partial
^{-}\left(  \frac{\partial u}{\partial x_{1}}\right)  ,v\right\rangle
-\left\langle v,\partial^{-}\left(  \frac{\partial u}{\partial x_{1}}\right)
\right\rangle +\int_{\left\{  \varepsilon\right\}  \times\Sigma}\left\langle
v,\partial^{-}u\right\rangle -\int_{\left\{  0\right\}  \times\Sigma
}\left\langle v,\partial^{-}u\right\rangle \\
&  =0.
\end{align*}

In order to estimate $\int_{A}\left\vert V_{x_{1}}\right\vert ^{2}$, we notice
that, for any fixed point $p\in\Sigma$, $v|_{\left[  0,\varepsilon\right]
\times\left\{  p\right\}  }$ can be treated as a function over the interval
$\left[  0,\varepsilon\right]  \ $and we compute,
\begin{align*}
\int_{0}^{\varepsilon}v^{2}dx_{1} &  =\int_{0}^{\varepsilon}\left(  \int
_{0}^{x_{1}}\frac{\partial v}{\partial x_{1}}\left(  t\right)  dt\right)
^{2}dx_{1}\\
&  \leq\int_{0}^{\varepsilon}\left(  \int_{0}^{x_{1}}ds\right)  \left(
\int_{0}^{x_{1}}\left\vert \frac{\partial v}{\partial x_{1}}\left(  t\right)
\right\vert ^{2}dt\right)  dx_{1}\\
&  \leq\int_{0}^{\varepsilon}x_{1}dx_{1}\int_{0}^{\varepsilon}\left\vert
\frac{\partial v}{\partial x_{1}}\left(  t\right)  \right\vert ^{2}dt\\
&  =\frac{\varepsilon^{2}}{2}\int_{0}^{\varepsilon}\left\vert \frac{\partial
v}{\partial x_{1}}\left(  t\right)  \right\vert ^{2}dt.
\end{align*}
This is basically an effective Poincar\'{e} inequality. By putting all these
together, we have
\begin{align*}
&  \left\langle \mathcal{D}V,\mathcal{D}V\right\rangle _{L^{2}}\\
&  =\int_{\left[  0,\varepsilon\right]  \times\Sigma}\left(  \left\Vert
u_{x_{1}}\right\Vert ^{2}+\left\Vert v_{x_{1}}\right\Vert ^{2}+\left\Vert
\partial^{-}v\right\Vert ^{2}+\left\Vert \partial^{+}u\right\Vert ^{2}\right)
\\
&  \geq\int_{0}^{\varepsilon}\int_{\Sigma}\left\Vert \partial^{+}u\right\Vert
^{2}+\int_{\Sigma}\int_{0}^{\varepsilon}\left\Vert v_{x_{1}}\right\Vert ^{2}\\
&  \geq\lambda_{\partial^{+}}\int_{0}^{\varepsilon}\int_{\Sigma}\left\Vert
u\right\Vert ^{2}+\frac{2}{\varepsilon^{2}}\int_{\Sigma}\int_{0}^{\varepsilon
}\left\Vert v\right\Vert ^{2}\\
&  \geq\min\left\{  \lambda_{\partial^{+}},2/\varepsilon^{2}\right\}  \left(
\int_{0}^{\varepsilon}\int_{\Sigma}\left\Vert u\right\Vert ^{2}+\int_{\Sigma
}\int_{0}^{\varepsilon}\left\Vert v\right\Vert ^{2}\right) \\
&  =\min\left\{  \lambda_{\partial^{+}},2/\varepsilon^{2}\right\}  \left\Vert
V\right\Vert _{L^{2}}^{2}.
\end{align*}
Therefore
\[
\lambda_{\mathcal{D}}\geq\min\left\{  \lambda_{\partial^{+}},2/\varepsilon
^{2}\right\}  .
\]

Conversely, we suppose $u$ is the first eigenfunction of $\Delta_{\Sigma
}=\partial^{-}\partial^{+}$, we extend $u$ trivially to a function on $\left[
0,\varepsilon\right]  \times\Sigma$ and define $V:=\left(  u,0\right)  \in
H_{-}^{1}$ then we have
\begin{align*}
\lambda_{\partial^{+}}\left\Vert V\right\Vert _{L^{2}}^{2} &  =\lambda
_{\partial^{+}}\int_{0}^{\varepsilon}\int_{\Sigma}\left\Vert u\right\Vert
^{2}\\
&  =\int_{0}^{\varepsilon}\int_{\Sigma}\left\Vert \partial^{+}u\right\Vert
^{2}\\
&  =\left\langle \mathcal{D}V,\mathcal{D}V\right\rangle _{L^{2}}\\
&  \geq\lambda_{\mathcal{D}}\left\Vert V\right\Vert _{L^{2}}^{2}.
\end{align*}
Combining these, we have the theorem.
\end{proof}

In order to adapt the simplified model above to our stituation, let us
consider $\mathcal{D}_{g}$ be the Dirac operator \ on a Riemannian Spin
manifold $\left(  A,g\right)  $ with metric $g.$ If we change the metric
conformally $g\rightarrow hg$ by any positive function $h\in C^{\infty}\left(
A\right)  $ then we have
\[
\mathcal{D}_{hg}=h^{-\frac{n+1}{4}}\circ\mathcal{D}_{g}\circ h^{\frac{n-1}{4}%
},
\]
where $n$ is the dimension of $A$. If we compare the Rayleigh quotient we
find
\[
\frac{1}{C}\inf_{V\in\mathbb{S}}\frac{\int_{X}\left\|  \mathcal{D}%
_{hg}V\right\|  _{hg}^{2}}{\int_{X}\left\|  V\right\|  _{hg}}\leq\inf
_{V\in\mathbb{S}}\frac{\int_{X}\left\|  \mathcal{D}_{g}V\right\|  _{g}^{2}%
}{\int_{X}\left\|  V\right\|  _{g}}\leq C\inf_{V\in\mathbb{S}}\frac{\int
_{X}\left\|  \mathcal{D}_{hg}V\right\|  _{hg}^{2}}{\int_{X}\left\|  V\right\|
_{hg}},
\]
where $C>0$ is a constant depending only on \ $0<\min_{x\in A}h\left(
x\right)  \leq\max_{x\in A}h\left(  x\right)  <\infty.$ \ In particular, this
implies%
\[
\frac{1}{C}\lambda_{\mathcal{D}_{g}}\leq\lambda_{\mathcal{D}_{hg}}\leq
C\lambda_{\mathcal{D}_{g}}.
\]

Now we can extend the above Theorem to any product three manifold
$A_{\varepsilon}=\left[  0,\varepsilon\right]  \times\Sigma$ with a
\textit{warped} product metric%
\[
g_{A_{\varepsilon}}=h\left(  x\right)  dx_{1}^{2}+g_{\Sigma}.
\]
This is because $g_{A_{\varepsilon}}$ is conformally equivalent to a product
metric $dx_{1}^{2}+h^{-1}g_{\Sigma}$ with conformal factor $h\left(  x\right)
$. Therefore we have the following corollary.

\begin{corollary}
\label{1Cor Eigen}Suppose $A_{\varepsilon}=\left[  0,\varepsilon\right]
\times\Sigma$ is equipped with a Riemannian metric of the form
$g_{A_{\varepsilon}}=h\left(  x\right)  dx_{1}^{2}+g_{\Sigma}$ for some smooth
positive function $h$ on $\Sigma.$Then we have
\[
\frac{1}{C}\lambda_{\partial^{+}}\leq\lambda_{\mathcal{D}}\leq C\lambda
_{\partial^{+}}%
\]
with some constant $C$ depend on $h.$
\end{corollary}

In particular we have the following corollary.

\begin{corollary}
Assumptions as before, we have $\ker\partial^{+}=0$ if and only if
$\ker\mathcal{D}=0.$
\end{corollary}

Remark: The above theorem gives an effective lower bound of the first
eigenvalue of $\mathcal{D}^{\ast}\mathcal{D}$. If we only need to obtain the
corollary we can also achieve this by exploring the complexified version of
the Cauchy-Riemann equations\ as follows
\begin{align*}
\frac{\partial}{\partial x_{1}}\left(  \frac{\partial v}{\partial x_{1}%
}\mathbf{i}\right)   &  =\frac{\partial}{\partial x_{1}}\left(  \bar{\partial
}_{z}u\right) \\
&  =\bar{\partial}_{z}\left(  \frac{\partial}{\partial x_{1}}u\right)
=\bar{\partial}_{z}\left(  -i\partial_{z}v\right) \\
&  =-\mathbf{i}\bigtriangleup v,
\end{align*}
and
\begin{align*}
\frac{\partial}{\partial x_{1}}\left(  \frac{\partial u}{\partial x_{1}%
}\mathbf{i}\right)   &  =\frac{\partial}{\partial x_{1}}\left(  \partial
_{z}v\right) \\
&  =\partial_{z}\left(  \frac{\partial}{\partial x_{1}}v\right)  =\partial
_{z}\left(  -\mathbf{i}\bar{\partial}_{z}u\right) \\
&  =-\mathbf{i}\bigtriangleup u\text{.}%
\end{align*}
Therefore we have
\[
\left(  \frac{\partial^{2}}{\partial x_{1}^{2}}+\bigtriangleup^{\Sigma
}\right)  u=0;\text{ \ }\left(  \frac{\partial^{2}}{\partial x_{1}^{2}%
}+\bigtriangleup^{\Sigma}\right)  v=0.
\]
We perform the Fourier transformation on $\Sigma,$ and we obtain
\[
\left(  \frac{\partial^{2}}{\partial x_{1}^{2}}-|\xi|^{2}\right)  \widehat
{v}=0
\]
subjected to the boundary conditions, $\widehat{v}|_{x_{1}=0}=\widehat
{v}|_{x_{1}=\varepsilon}=0$. Clearly\ this implies $v=0$ and therefore
\[
\frac{\partial u}{\partial x_{1}}=\bar{\partial}u=0.
\]
That is $u\in\ker\bar{\partial}_{z.}$

\subsection{Perturbation arguments}

Let $C\subset M$ be a coassociative submanifold. Suppose that $v$ is a normal
vector field on $C$ such that its corresponding self-dual two form, $\eta
_{0}=\iota_{v}\Omega\in\wedge_{+}^{2}\left(  C\right)  $ \ is harmonic with
respect to the induced metric. So $\eta_{0}$ is actually a symplectic form on
the complement of the zero set $Z\left(  \eta_{0}\right)  $ of $\eta_{0}$ in
$C$. Furthermore\ $v$ defines an almost complex structure $J_{v}$ on the
$C\backslash Z\left(  \eta_{0}\right)  $ via $\left|  v\right|  ^{-1}%
v\times\cdot.$ Since the deformation of coassociative manifold is
unobstructed, we may assume that there is an one parameter family of
coassociative submanifolds $C_{t}$ in $M$ which corresponds to integrating out
the normal vector field $v,$ that is
\[
\left.  \frac{dC_{t}}{dt}\right|  _{t=0}=v\in\Gamma\left(  C_{0},N_{C_{0}%
/M}\right)  .
\]

In this article we assume that $\eta_{0}$ is nowhere vanishing on $C$, that is
$\left(  C,\eta_{0}\right)  $\ is a symplectic four manifold. We are going to
establish a relation between the non-vanishing of Seiberg-Witten invariants of
$C$ and the existence of instantons with coassociative boundary conditions.

Given any submanifold $\Sigma$ in $C$, using the variation normal vector
fields $dC_{t}/dt$ to identify various $C_{t}$'s, we obtain submanifolds
$\Sigma_{t}$ in $C_{t}$ vary smoothly with respect to $t$. We denote
\[
A_{\varepsilon}^{\prime}:=\bigcup_{0\leq t\leq\varepsilon}\Sigma_{t}%
\]
Since $\eta$ is nonvanishing on $C$, all coassociative submanifolds$\ C_{t}$'s
are mutually disjoint. In particular $A_{\varepsilon}^{\prime}$ is a smooth
three dimensional submanifold in $M$. Suppose that $\Sigma\subset C$ is a
$J_{v}$-holomorphic curve, then the tangent spaces of $A_{\varepsilon}%
^{\prime}$ is associative along $\Sigma$. In fact $A_{\varepsilon}^{\prime}$
is close to be an instanton when $\varepsilon$ is small, more precisely we
have the following:

(i) $\left|  \tau|_{A_{\varepsilon}^{\prime}}\right|  \leq c_{1}\varepsilon
,$for some constant $c_{1}$.

(ii) The natural diffeomorphism between $A_{\varepsilon}^{\prime}$ and
$\Sigma\times\left[  0,\varepsilon\right]  $ is a $\varepsilon$-isometry
between the induced metric $g_{A_{\varepsilon}^{\prime}}$ on $A_{\varepsilon
}^{\prime}$ and the \textit{warped} product metric
\[
g_{\Sigma}+h\left(  x\right)  dt^{2}%
\]
on $\Sigma\times\left[  0,\varepsilon\right]  $, where $g_{\Sigma}$ is the
induced metric on $\Sigma$ and $h\left(  x\right)  $ is the length of
$v=dC_{t}/dt|_{t=0}$ restricted to $\Sigma$.

(iii) The derivative $F^{\prime}\left(  0\right)  $ of the functional
\[%
\begin{array}
[c]{ccc}%
F:\Gamma\left(  A_{\varepsilon},N_{A_{\varepsilon}/M}\right)  &
\longrightarrow & \Gamma\left(  A_{\varepsilon},N_{A_{\varepsilon}/M}\right)
\\
V & \longrightarrow & \ast\left(  \exp_{V}^{\ast}\tau\right)  ^{\perp}%
\end{array}
\]
is $\varepsilon$-close to the twisted Dirac operator on $A_{\varepsilon}$ with
respect to the above warped product metric. That is,
\[
F^{\prime}\left(  0\right)  :H^{1}\left(  N_{A_{\varepsilon}^{\prime}%
/M}\right)  \rightarrow L^{2}\left(  N_{A_{\varepsilon}^{\prime}/M}\right)
\text{ and }\left|  F^{\prime}\left(  0\right)  -\mathcal{D}_{\Sigma
\times\left[  0,\varepsilon\right]  }\right|  _{C^{0}}\leq c_{2}%
\varepsilon\text{,}%
\]
for some constant $c_{2}$. Note that we need to use the orthogonal projection
to identify the normal bundle to $A_{\varepsilon}^{\prime}$ to the twisted
spinor bundle of $\Sigma\times\left[  0,\varepsilon\right]  $.

These imply that $\left|  \lambda\left(  F^{\prime}\left(  0\right)  \right)
-\lambda\left(  \mathcal{D}\right)  \right|  \leq\frac{1}{2}\lambda\left(
\mathcal{D}\right)  $ for $\varepsilon<\delta$ where $\delta$ is some small
positive number depend on the geometry of $C_{0}$. Recall from proposition
\ref{1Prop alm instanton} that $F\left(  V\right)  =0$ implies that the image
of the exponential map $A=\exp_{V}\left(  A_{\varepsilon}\right)  $ is an
instanton in $M$. To find the zeros of $F$, we are going to apply the
following quantitative version of the implicit function theorem.

\begin{theorem}
Let $X$ and $Y$ be Banach space and $F:B_{r}\left(  x_{0}\right)  \subset
X\rightarrow Y$ a $C^{1}$-map, such that

\begin{enumerate}
\item $\left(  DF\left(  x_{0}\right)  \right)  ^{-1}$ is a bounded linear
operator with $\left|  \left(  DF\left(  x_{0}\right)  \right)  ^{-1}F\left(
x_{0}\right)  \right|  \leq\alpha$ and $\left|  \left(  DF\left(
x_{0}\right)  \right)  ^{-1}\right|  \leq\beta;$

\item $\left|  DF\left(  x_{1}\right)  -DF\left(  x_{2}\right)  \right|  \leq
k\left|  x_{1}-x_{2}\right|  $ for all $x_{1},x_{2}\in B_{r}\left(
x_{0}\right)  ;$

\item $2k\alpha\beta<1$ and $2\alpha<r.$

Then $F$ has a unique zero $z$ in $B_{r}\left(  x_{0}\right)  .$
\end{enumerate}
\end{theorem}

By combining with the eigenvalue estimates in section \ref{1Sec e v estimate},
we can obtain the following result on the existence of instantons.

\begin{theorem}
Suppose that $M$ is a $G_{2}$-manifold and $C_{t}$ is an one parameter family
of coassociative submanifolds in $M$ with $\left(  dC_{t}/dt\right)  |_{t=0}$ nonvanishing.

For any regular $J$-holomorphic curve $\Sigma$ in $C_{0}$, there is an
instanton $A_{\varepsilon}$ in $M$ which is diffeomorphic to $\left[
0,1\right]  \times\Sigma$ and $\partial A_{\varepsilon}\subset C_{0}\cup
C_{\varepsilon}$, for all sufficiently small positive $\varepsilon$.
\end{theorem}

\begin{proof}
Consider the functional $F_{\varepsilon}:X\rightarrow Y$ with $X=H^{1}\left(
N_{A_{\varepsilon}^{\prime}/M}\right)  $, $Y=L^{2}\left(  N_{A_{\varepsilon
}^{\prime}/M}\right)  $ and $F_{\varepsilon}\left(  V\right)  =\ast\left(
\exp_{V}^{\ast}\tau\right)  ^{\perp}$. In order to apply the above implicit
function theorem, we need to know that
\[
\left|  \left(  DF_{\varepsilon}\left(  0\right)  \right)  ^{-1}%
F_{\varepsilon}\left(  0\right)  \right|  \leq\alpha\text{.}%
\]
From Theorem \ref{1Thm Eigen} and its Corollary \ref{1Cor Eigen} we know that
for $\varepsilon$ small we have
\[
\left|  \left(  DF_{\varepsilon}\left(  0\right)  \right)  ^{-1}\right|  \leq
C\sqrt{\frac{2}{\lambda_{\bar{\partial}}}}%
\]
Also by our construction we have
\[
\lim_{\varepsilon\rightarrow0}F_{\varepsilon}\left(  0\right)  =0
\]
so for $\varepsilon$ small, we have \
\[
\lim_{\varepsilon\rightarrow0}\left|  \left(  DF_{\varepsilon}\left(
0\right)  \right)  ^{-1}F_{\varepsilon}\left(  0\right)  \right|  =0.
\]
By applying the above implicit function theorem and proposition \ref{1Prop alm
instanton} we obtain an one parameter family of instantons $A_{\varepsilon}$
for $\varepsilon$ small, and \ $\partial A_{\varepsilon}\subset C_{0}\cup
C_{\varepsilon}.$
\end{proof}

\bigskip

In particular, by combining Taubes' results with the above theorem we obtain
the following existence result.

\begin{corollary}
Suppose that $C$ is a coassociative submanifold in a $G_{2}$-manifold $M$ with
non-trivial Seiberg-Witten invariants. Given any symplectic form on $C$, we
write $C_{t}$'s the corresponding coassociative deformations of $C$ in $M$.
Then there is an instanton $A_{t}$ in $M$ with boundaries lying on $C_{0}\cup
C_{t}$ for each sufficiently small $t$.
\end{corollary}

Lastly we expect that any instanton $A$ in $M$ bounding $C_{0}\cup C_{t}$ and
with small volume must arise in the above manner. Namely we need to prove a
$\varepsilon$-regularity result for instantons.

\bigskip

\textit{Acknowledgments: The first author is partially supported by
NSF/DMS-0103355 and he expresses his gratitude to J.H. Lee, Y.G. Oh, C.
Taubes, R. Thomas and A. Voronov for useful discussions. The second author
thanks S.L. Kong, G. Liu, Y.J. Lee for useful discussions. }

\bigskip

{\small Addresses:}

{\small N.C. Leung}

{\small School of Mathematics, University of Minnesota, Minneapolis, MN 55454}

{\small Email: LEUNG@MATH.UMN.EDU}

{\small \bigskip}

{\small X.W. Wang}

{\small Department of Mathematics, University of California, Los Angeles, CA 90095}

{\small Email: XIAOWEI@MATH.UCLA.EDU}
\end{document}